\documentclass[graybox]{svmult}

\usepackage{mathptmx}       
\usepackage{helvet}         
\usepackage{courier}        
\usepackage{type1cm}        
\usepackage{makeidx}         
\usepackage{graphicx}        
\usepackage{multicol}        
\usepackage[bottom]{footmisc}

\makeindex             

\newcommand{\s}{\mathcal}

\def \-{\bar}
\newcommand{\p}{\partial}

\date{}

\begin{document}

\title*{\bf On the CR transversality of holomorphic  maps into hyperquadrics}

\bigskip
\subtitle{\small Dedicated to Professor Yum-Tong Siu on the
occasion of his 70th birthday}

\author{Xiaojun Huang \ and \ Yuan Zhang}
\institute{Xiaojun Huang, partially supported in part by National
Science Foundation DMS-1363418.  \at Department of Mathematics,
Rutgers University, New Brunswick, NJ 08903, USA
\email{huangx@math.rutgers.edu} \and Yuan Zhang, partially supported
in part by National Science Foundation DMS-1265330. \at Department
of Mathematical Sciences, Indiana University - Purdue University,
Fort Wayne, IN 46805, USA \email{zhangyu@ipfw.edu}} \maketitle

\abstract{ Let $M_\ell$ be  a smooth Levi-nondegenerate hypersurface
of signature $\ell$ ($0<\ell<\frac{n-1}{2}$) in $\mathbf C^n$ with $
n\ge 3$, and write $H_\ell^N$ for the standard hyperquadric of the
same signature in $\mathbf C^N$ with $N-n< \frac{n-1}{2}$. Let $F$
be a holomorphic map sending $M_\ell$ into $H_\ell^N$. Assume $F$
does not send a neighborhood of $M_\ell$ in $\mathbf C^n$ into
$H_\ell^N$.  We show that $F$  is necessarily CR transversal to
$M_\ell$ at any point. Equivalently, we show that $F$ is a local CR
embedding from $M_\ell$ into $H_\ell^N$.}

\bigskip
\section{Introduction  and the main theorems}
Let $M_1$ and $M_2$ be two connected   smooth CR hypersurfaces in
$\mathbf C^n$ and $\mathbf C^N$, respectively, with $3\le n\le N$.
Let $F$ be a holomorphic map from some small neighborhood $U\subset
\mathbf C^{n}$ of $M_1$ into $\mathbf C^N$ with $F(M_1)\subset M_2$.
Given a point $p\in M_1$, denote by $T^{(1, 0)}_p M_1$ and $T^{(1, 0)}_{F(p)} M_2$ the holomorphic
tangent vector spaces of $M_1$ at $p$ and $M_2$ at $F(p)$, respectively. Assume $F$ does not send a
neighborhood of $p$ in $\mathbf C^n$ into $M_2$. An important
question in the study of the geometric structure of $F$   is to
understand the geometric conditions for the manifolds in which $F$
is CR transversal to $M_1$ at $p$. Recall that  $F$ is said to be CR
transversal at $p$ if
$$
 T^{(1,0)}_{F(p)}M_2+dF(T^{(1,0)}_p\mathbf C^{n})=T^{(1,0)}_{F(p)}\mathbf C^{N}.
$$
Roughly speaking, the CR transversality property can be interpreted
as an non-vanishing property of  the normal derivative of the normal
components for the map.

The problem has been extensively investigated in the literature.
When both the target and the source manifolds are strongly
pseudoconvex, CR transversality always holds due to the classical
Hopf lemma. In the equal dimensional case ($n=N$), work has been
done by Pinchuk \cite{Pi}, Fornaess \cite{Fo}, Baouendi-Rothschild
\cite{BR}, Ebenfelt-Rothschild \cite{ER}, Huang \cite{Hu2}, Isaev
[Iz1] [Iz2], Huang-Pan [HP] and the references therein. The study of
the higher codimensional case starts with the work of Baouendi-Huang
in \cite{BH} where it is proved that the CR transversality always
holds when the manifolds are hyperquadrics of the same signature.
Baouendi-Ebenfelt-Rothschild \cite{BER2} proved, under rather
general setting,  that the CR transversality holds in an open dense
subset. See also a recent paper of Ebenfelt-Son \cite{ES} and the
references therein.

While there exist examples where CR transversality fails on certain
thin sets (see, for instance \cite{BER2}), as mentioned above,  the
rigidity theorem  due to Baouendi-Huang \cite{BH} indicates that the
CR transversality holds everywhere when both $M_1$ and $M_2$ are
hyperquadrics of the same signature $\ell$. Enlightened by this
result,  the following conjecture concerning the CR transversality
was asked by Baouendi and the first author in the year of 2005:

\medskip
{\bf Conjecture (Baounedi-Huang, 2005)}: {\it Let $M_1\subset
{\mathbf C}^{n} $ and $M_2\subset {\mathbf C}^{N}$ be two
(connected) Levi non-degenerate real analytic hypersurfaces with the
same signature $0<\ell\le\frac{n-1}{2}$. Here $3\le n<N$. Let $F$ be
a holomorphic map defined in a neighborhood $U$ of $M_1$, sending
$M_1$ into $M_2$. Then either $F$ is a local CR embedding from $M_1$
into $M_2$ or $F$ is totally degenerate in the sense that it maps  a
neighborhood $U$ of $M_1$ in ${\mathbf C}^{n}$ into $M_2$.}

\medskip

We point out  that, for the $M_1$ and $M_2$ given in the conjecture,
the fact that $F$ is CR transversal at $p$ is equivalent to the fact
that $F$ is  a CR embedding from a neighborhood of $p$ in $M_1$ into
$M_2$. Along these lines, in a recent paper of the  authors
\cite{HZ2}, by developing a new technique, we showed the CR
transversality holds when  $M_2 =H_\ell^{n+1}$ and the point under
study is not CR umbilical in the sense of Chern-Moser. Recall that a hyperquadric $H_\ell^n$ of signature $\ell$ in $\mathbf C^n$ is defined by
$$
  H_\ell^n:= \big\{(z, w)\in\mathbf C^{n-1}\times\mathbf C: \Im w =|z|^2_\ell\big\},
$$
where for any $n$-tuples $a$ and $b$, $\langle a, \bar b\rangle_{\ell}: = -\sum\limits_{j=1}^\ell a_j\bar b_j+ \sum\limits_{j=\ell+1}^n a_j\bar b_j$ and $|a|^2_\ell = \langle a, \bar a\rangle_{\ell}$.

In this paper, combining  a quantitative version of a very useful
lemma due to the first author with the tools developed in
\cite{HZ2}, we are able to drop the geometric assumption of the
umbilicality and relax the codimension-one restriction in
\cite{HZ2}. The generalization of the above mentioned  lemma in
\cite{Hu} will be addressed in detail in section 3.

We next state our main theorems:
\medskip

\begin{theorem}\label{main}
 Let $M_{\ell}$ be a smooth Levi non-degenerate hypersurface of signature $\ell$ in $\mathbf C^{n}$ with $n\ge 3$ and $0\in
 M_\ell$. (We assume $0<\ell\le \frac{n-1}{2}$.)
 Suppose that $F$ is a holomorphic map  in a small neighborhood
    $U$ of $0\in {\mathbf
   C}^{n}$ such that
    $$F(M_{\ell}\cap U)\subset H_{\ell}^{N}$$ with $N-n< \frac{n-1}{2}$. If  $F(U)\not \subset H_{\ell}^{N}$,
    then $F$ is  CR transversal to $M_\ell$ at $0$, or equivalently,
$F$ is a CR embedding from a  small neighborhood of $0\in M_\ell$
into $H_{\ell}^{N}$.
\end{theorem}
\medskip

\begin{theorem}\label{21}
Let $M_{\ell}$ be a germ of a smooth Levi non-degenerate
hypersurface at 0 of signature $\ell$ in $\mathbf C^{n}$, $n\ge 3$.
(We assume $0<\ell\le \frac{n-1}{2}$.) Suppose that there exists   a
holomorphic map $F$ in a neighborhood $U$ of $0$ in $\mathbf C^{n}$
sending $M_{\ell}$ into $H_{\ell}^{N}$ but $F(U)\not\subset
H_{\ell}^{N}$, $N< 2n-1$. Then $M_{\ell}$ is CR embeddable into
$H_{\ell}^{N}$ near $0$. Equivalently,  there exists a holomorphic
map $\tilde F: M_\ell \rightarrow H_{\ell}^{N}$ near $0$, which is
CR transversal to $M_\ell$ at $0$.
\end{theorem}
\medskip

The idea of the proof is based on a re-scaling technique that was
initially introduced in \cite{HZ2}. With the aid of a quantitative
lemma of the first author in \cite{Hu}, we  generate a formal CR
transversal map which, by a result of Meylan-Mir-Zaitsev proved in
\cite{MMZ}, is necessarily convergent. Finally, using a rigidity
result in \cite{EHZ}, $F$
 differs from the CR transversal map only by an automorphisms of
the target and hence it is CR transversal as well.

The outline of the paper is as follows. In section 2, the notations
and a normalization procedure of Baouendi-Huang is revisited. A
modified  lemma in \cite{Hu} is discussed and proved in section 3.
Section 4 is devoted to the proof of the main theorem.

\medskip

\section{Notations and a normalization procedure}

Let $M_\ell$ be a germ at 0 of a smooth Levi non-degenerate hypersurface of signature $\ell$ in $\mathbf C^{n}$. After a holomorphic change of coordinates, $M_\ell$ near the origin can be expressed as follows.
 \begin{equation}\label{MM}
    M_\ell =\big\{(z, w)\in\mathbf C^{n-1}\times\mathbf C: \Im w =|z|^2_\ell-\frac{1}{4}\s S(z)+o(4)\big\}.
 \end{equation}
  Here $\s S(z): = \sum\limits_{1\le\alpha,{\beta}, \gamma, \delta\le n} s_{\alpha\bar {\beta}\gamma \bar{\delta}}z_{\alpha} {\bar z_\beta}z_{\gamma}{\bar
z_\delta}$ is a  homogeneous polynomial of bi-degree (2,2), called
the Chern-Moser-Weyl curvature function of $M_\ell$ at 0. See
\cite{CM} for more details. In what follows, we always assume that
$\ell\le (n-1)/2$. Hence $\ell$ is a holomorphic  invariant.
\medskip

As in \cite{CM}, assign the weighted degree 1 to variable $z$  and 2 to variable $w$. Given a holomorphic function $h$, denote by $h^{(k)}$ the terms of weighted degree $k$, and by $h^{(\mu, \nu)}$ the terms of degree $\mu$ in $z$ variable and of degree $\nu$ in $w$ variable in the power series expansion of $h$ at 0.
For each integer $k\ge 0$, we  write $o(k)$ for terms of degree larger than $k$, and $o_{wt}(k)$ for terms of weighted degree larger than $k$.  To
simplify our notation, we also preassign the coefficient of $h$ with
negative degrees to be 0.
 \medskip

Let $\tilde M_\ell$ be a germ at 0 of another smooth Levi-nondegenerate hypersurface in $\mathbf C^{N}$ of  signature $\ell$ given by
\begin{equation}\label{MMM}
    \tilde M_\ell =\big\{(\tilde z, \tilde w)\in\mathbf C^{N-1}\times\mathbf C: \Im \tilde w =|\tilde z|^2_\ell-\frac{1}{4} \tilde {\s S}(\tilde z)+o(4)\big\}.
 \end{equation}
 Here $\tilde {\s S}$ is the corresponding Chern-Moser curvature tensor function of $\tilde M_\ell$ at $0$.

  Let $F$ be a smooth CR map sending $(M_\ell, 0)$ into $(\tilde M_\ell, 0)$. Write
 \begin{equation}\label{F}
 F: = (\tilde f, g)=(f, \phi, g)
  \end{equation}
with $f=({f_1}, \ldots, {f_{n-1}})$ and $\phi = ( \phi_1, \ldots, \phi_{N-n})$ being components of $F$.  Assume that $F$ is CR transversal at $0$. Then, following a normalization procedure as in [$\S 2$, BH],
we have
\begin{eqnarray}\label{eqn:022}
&\tilde z=(f_1(z, w), \ldots, f_{n-1}(z, w), \phi_1(z, w), \ldots, \phi_{N-n}(z, w))=\lambda
z U+\vec{a}w+ O(|(z, w)|^2)\nonumber\\
&\tilde w=g(z, w)=\sigma\lambda^2w+O(|(z, w)|^2).
\end{eqnarray}
Here $U$ can be extended to an ${(N-1)}\times {(N-1)}$ matrix
$\tilde U\in SU(N-1, \ell)$ (namely $\langle X\tilde U,
Y\overline{\tilde U}\rangle_\ell=\langle X, Y\rangle_\ell$ for any $X, Y\in
\mathbf C^{N-1}$), $\ \vec{a}\in \mathbf C^{N-1}$ and $\lambda
>0$, $\sigma=\pm1$ with $\sigma= 1$ for $\
\ell<\frac{n-1}{2}$. When $\sigma=-1$, by considering
$F\circ\tau_{{n-1}/2}$ instead of $F$, where
 $\tau_{\frac{n-1}{2}}(z_1,\ldots,
z_{\frac{n-1}{2}},z_{\frac{n-1}{2}+1},\ldots,
 z_{n-1}, w)=(z_{\frac{n-1}{2}+1},\ldots,
 z_{n-1},z_1,\ldots,
z_{\frac{n-1}{2}},-w),$ we can make $\sigma=1$.
Hence, we will assume
in what follows that $\sigma=1$. Moreover, as in  [HZ],
$F$ can be  normalized as follows:
\medskip

\begin{proposition} ({\cite{HZ}})\label{55-} Let $M_\ell$ and $\tilde M_\ell$ be defined by
(\ref{MM}) and (\ref{MMM}), respectively, and let $F$ be a smooth CR map sending $M_\ell$ into
$\tilde M_\ell$ given by (\ref{F}) and (\ref{eqn:022}) with $\sigma=1$. Then after composing $F$  from the left by some automorphism $T\in Aut_0(H_\ell^{N})$ preserving the origin, the following holds:

 $$
     F^\sharp =(f^\sharp, \phi^\sharp, g^\sharp): =T\circ F,
  $$
   with
    \begin{eqnarray*}
       &{f^\sharp} (z, w)=z+\frac{i}{2}a^{(1,0)} (z)w +o_{wt}(3),\\
    &\phi^\sharp (z, w)=\phi^{(2,0)} (z) + o_{wt}(2),\\
    &g^\sharp (z, w)=w+o_{wt}(4),
\end{eqnarray*}
and
$$
    \langle a^{(1,0)}(z), \bar z\rangle_{\ell} |z|_{\ell}^2 = |\phi^{(2,0)}(z)|^2+ \frac{1}{4}(\s S(z)-\lambda^{-2}\tilde{\s S}(\lambda(z, 0)\widetilde {U})).
$$
In particular,   the automorphism $T$ is given by
$$
T(\tilde z, \tilde w)=\frac{(\lambda^{-1}(\tilde z- \lambda^{-2}\vec
a \tilde w)\tilde U^{-1}, \lambda^{-2}\tilde w)}{q(\tilde z,
\tilde w)}
$$
with
$r_0=\frac{1}{2}\Re \{g^{''}_{ww}(0)\},\ q(\tilde z, \tilde
w)=1+2i\langle\tilde z,  \lambda^{-2}\overline{\vec
a}\rangle_\ell+\lambda^{-4}(r_0-i|\vec a|_\ell^2)\tilde w$.
Moreover, $F^\sharp$ sends $M_\ell$ into $\tilde M^\sharp: =T(\tilde M_\ell)$ given by
$$
 \tilde M^\sharp = \{(\tilde z^\sharp, \tilde w^\sharp)\in\mathbf C^{N+1}: \Im \tilde w^\sharp = |\tilde z^\sharp|_\ell^2+
  \frac{1}{4}\tilde{\s S}^\sharp(\tilde z^\sharp)+o(4) \big\}
$$
with $\tilde {\s S}^{\sharp}(z^{\sharp})= \lambda^{-2}\tilde {\s
S}(\lambda z^{\sharp}\tilde U)$.
\end{proposition}
\medskip

\section{A quantitative  version of a basic lemma}

In this section,  some simple preparation facts will  be given
without proof at first. In the second part of the section, we will
discuss a quantitative version of a  lemma obtained in \cite{Hu},
which played crucial role for us to get the convergence in our
rescaling argument.

Given a  polynomial $\phi$, define $\|\phi\|$ to be the maximum modulus of all the coefficients in $\phi$. For a given vector-valued polynomial $\phi=(\phi_1, \ldots, \phi_s)$, $\|\phi\|: =\max_{1\le j\le s}\|\phi_j\|$. We first refer to a lemma in \cite{HZ2} without proof.
\medskip

\begin{lemma}\label{18}\cite{HZ2}
 (1). Let $X(z,\-{z})$ and $Y(z,\-{z})$ be two polynomials
 such that $X(z,\-{z})=Y(z,\-{z})|z|_\ell^2.$ Then $\|Y\|$ is bounded by a constant depending only on  $\|X\|$ and the degree of $X$.

 (2). Let $h(z)$ be a homogeneous holomorphic polynomial of degree $d$ in $z\in\mathbf C^n$. If $|h(z)|\le  c |z|^d$ on
   $\{|z|_{\ell}^2=0\}$, then $\| h\|\le C$ for some $C$ depending only  on $c$ and $d$.
\end{lemma}
\medskip

In various rigidity problems concerning CR immersions, the following
lemma in \cite{Hu} plays an essential role in deriving key
identities to eventually conclude uniqueness:
\medskip

  \begin{lemma}\label{Hu}\cite{Hu}
  Let $\{\phi_j\}_{j=1}^{n-1}$ and $\{\psi_j\}_{j=1}^{n-1}$ be two families of holomorphic functions in $\mathbf C^n$. Let $B(z, \xi)$ be a real-analytic function in $(z, \xi)$. Suppose that
$$
  \sum\limits_{j=1}^{n-1}\phi_j(z)\psi_j(\xi)= B(z, \xi)\langle z, \xi\rangle_\ell.
$$
Then $B(z, \xi) = \sum\limits_{j=1}^{n-1}\phi_j(z)\psi_j(\xi) =0$.
  \end{lemma}
\medskip

We find a quantitative version of the above lemma serves our purpose
under this context perfectly well.
\medskip

\begin{lemma}\label{bound}
  Let $\{\phi_j\}_{j=1}^{n-1}$ and $\{\psi_j\}_{j=1}^{n-1}$ be two families of holomorphic polynomials of degree $k$ and $m$ in $\mathbf C^n$, respectively. Let $H(z, \xi), B(z, \xi)$ be two polynomials in $(z, \xi)$. Suppose that
$$
  \sum\limits_{j=1}^{n-1}\phi_j(z)\psi_j(\xi)= H(z, \xi) + B(z, \xi)\langle z, \xi\rangle_\ell
$$
and $\|H\|\le C$. Then $\| B\| \le \tilde C$ and $
\|\sum\limits_{j=1}^{n-1}\phi_j(z)\psi_j(\xi)\|\le \tilde C$ with $\tilde
C$ dependent only on $(C, k, m, n)$.
\end{lemma}
\medskip

The proof of the lemma is based on the following algorithm together
with Lemma \ref{Hu}.  First, let us formulate the algorithm
procedure so as to re-adjust two families $\{\phi_j\}_{j=1}^{n-1}$
and $\{\psi_j\}_{j=1}^{n-1}$ in Lemma \ref{bound}.
\medskip

\begin{lemma}\label{6}
Let $\phi: =\{\phi_j\}_{j=1}^{s}$ and $\psi: = \{\psi_j\}_{j=1}^{s}$ be two families of holomorphic polynomials of degree $k$ and $m$ in $\mathbf C^n$, respectively.
 There exist  two families $\tilde\phi: =\{\tilde\phi_j\}_{j=1}^{s}$ and $\tilde \psi: =\{\tilde \psi_j\}_{j=1}^{s}$ of holomorphic polynomials of degree $k$ and $m$ in $\mathbf C^n$, respectively, such that
\begin{equation}\label{16}
  \sum\limits_{j=1}^{s}\phi_j(z)\psi_j(\xi)=\sum\limits_{j=1}^{s}\tilde\phi_j(z)\tilde\psi_j(\xi)
\end{equation}
and
\begin{equation}\label{17}
 \|\tilde\phi\|\le 1, \ \ C \|\tilde\psi\|\le \|\sum\limits_{j=1}^{s}\tilde\phi_j(z)\tilde\psi_j(\xi)\|\le  s\|\tilde\psi\|
\end{equation}
for some positive constant $C$ dependent only on $s$.
\end{lemma}
\medskip

\noindent{\bf Proof of Lemma \ref{6}:} Without loss of generality, assume $\|\phi_j\|\ne 0$ for all $1\le j\le s$ and $\{\phi_j\}_{j=1}^s$ are linearly independent. Moreover, by replacing $\phi_j$ and $\psi_j$ by $\frac{\phi_j}{\|\phi_j\|}$ and $\|\phi_j\|\psi_j$, respectively, one can  assume that $\|\phi_j\|=1$ for all $1\le j\le s$.
Denote by $\{e_l\}_{l=1}^{d(k)}$  a basis of unit monomials  to span the polynomial spaces of degree $k$ and write $\phi_j=\sum\limits_{1\le l\le d(k)}D_j^l e_l, 1\le j \le s$. Here $d(k)$ is the dimension of polynomial spaces of degree $k$. Hence $\|\phi_j\|=\max\limits_{1\le l\le d(k)} D_j^l$ for each $1\le j\le s$. Arranging the order of $\{e_l\}$ if necessary, we can make $D_1^1= 1$ and $|D_1^l|\le 1$.

{\bf Step 1:} Let $^1\phi_1: =\phi_1, ^1\phi_j: =\phi_j-D_j^1\cdot\phi_1, 2\le j\le s$. Then in terms of the basis representation $^1\phi_j: =^1D_j^l\cdot e_l$, one has
\begin{eqnarray*}
    &^1D_1^1=1, \ |^1D_1^l|\le 1, \ 2\le l\le d(k);\\
  &^1D_j^1=0, \ |^1D_j^l|\le 2, \ 2\le j\le s, \ \ 2\le l\le d(k).
\end{eqnarray*}
Moreover, letting $^1\psi_1: = \psi_1+\sum\limits_{j=2}^s D_j^1\cdot\psi_j, ^1\psi_j: = \psi_j, 2\le j \le s,$ then
\begin{equation}\label{101}
 \sum\limits_{j=1}^{s}\phi_j(z)\psi_j(\xi)=\sum\limits_{j=1}^{s}{^1\phi_j(z)}\cdot{^1\psi_j(\xi)}.
\end{equation}

{\bf Step 2:} Normalize $^1\phi_j, 2\le j \le s$ by replacing $^1\phi_j, ^1\psi_j$  by $\frac{^1\phi_j}{\|^1\phi_j\|}$ and $\|^1\phi_j\|\cdot{^1\psi_j}$, respectively.  By abuse of notation, we still denote them by $^1\phi_j, ^1\psi_j$ and the representation matrix under the basis $\{e_l\}$ by $\{^1D_j^l\}$. Moreover, since $\{\phi_j\}_{j=1}^s$ are linearly independent, by rearranging the order of $\{e_l\}_{l=2}^{d(k)}$ if necessary, we have (\ref{101}) holds with
\begin{eqnarray*}
     &^1D_1^1=1, \ |^1D_1^l|\le 1, \ 2\le l\le d(k);\\
     &^1D_2^1=0, \ ^1D_2^2=1, \ |^1D_2^l|\le 1, \ 3\le l\le d(k);\\
    &^1D_j^1=0, \ |^1D_j^l|\le 1, \ 3\le j\le s, \ \ 2\le l\le d(k)
\end{eqnarray*}
and for each $1\le j\le s$,
$$ \max\limits_{1\le l\le d(k)} {^1}D_j^l =1.$$

{\bf Step 3:} Let $^2\phi_2={^1}\phi_2, {^2}\phi_j: = {^1}\phi_j - {^1}D_j^2\cdot{^1}\phi_2$ for $1\le j\le s, j\ne 2$. Then in terms of the basis representation $^2\phi_j: ={^2}D_j^l\cdot e_l$, we deduce
\begin{eqnarray*}
    &^2D_1^1=1, \ ^2D_1^2=0,\  |^2D_1^l|\le 2, \ 3\le l\le d(k);\\
     &^2D_2^1=0, \ ^2D_2^2=1, \  |^2D_2^l|\le 1, \ 3\le l\le d(k);\\
    &^2D_j^1=0, \ ^2D_j^2=0,\  |^2D_j^l|\le 2, \ 3\le j\le s, \ \ 3\le l\le d(k).
\end{eqnarray*}
Moreover, letting $^2\psi_2: = {^1}\psi_2+\sum\limits_{j\ne 2} {^1}D_j^2\cdot{^1}\psi_j, {^2}\psi_j: = {^1}\psi_j, 1\le j \le s$ with $j\ne 2$, then
\begin{equation}\label{102}
 \sum\limits_{j=1}^{s}\phi_j(z)\psi_j(\xi)=\sum\limits_{j=1}^{s}{^2\phi_j(z)}\cdot{^2\psi_j(\xi)}.
\end{equation}

{\bf Step 4:} Normalize $^2\phi_j, 1\le j \le s, j\ne 2$ by replacing $^2\phi_j, ^2\psi_j$  by $\frac{^2\phi_j}{\|^2\phi_j\|}$ and $\|^2\phi_j\|\cdot{^2\psi_j}$, respectively. As before, we still denote them by $^2\phi_j, ^2\psi_j$ and the representation matrix under the basis $\{e_l\}$ by $\{^2D_j^l\}$. Furthermore, (\ref{102}) holds with
\begin{eqnarray*}
    &1\ge{^2}D_1^1\ge \frac{1}{2}, \ ^2D_1^2=0,\  |^2D_1^l|\le 1, \ 3\le l\le d(k);\\
     &^2D_2^1=0, \ ^2D_2^2=1, \  |^2D_2^l|\le 1, \ 3\le l\le d(k);\\
    &^2D_j^1=0, \ ^2D_j^2=0,\  |^2D_j^l|\le 1, \ 3\le j\le s, \ \ 3\le l\le d(k)
\end{eqnarray*}
and for each $1\le j\le s$,
$$ \max\limits_{1\le l\le d(k)} {^2}D_j^l =1.$$

{\bf Step 5:} Continue the above process until we get new families $\{^s\phi_j\}_{j=1}^s, \{^s\psi_j\}_{j=1}^s$ such that under the basis representation, $^s\phi_j: ={^s}D_j^l\cdot e_l$ with
\begin{eqnarray*}
 ^sD=\left[\matrix{
   {^s}D_1^1 &0&0       &\cdots &0  &^sD_1^{s+1} &\cdots &^sD_1^{d(k)}\\
  0       &^sD_2^2 &0       &\cdots &0   &^sD_2^{s+1} &\cdots&^sD_2^{d(k)}\\
  0       &0       &^sD_3^3 &\cdots &0   &^sD_2^{s+1} &\cdots &^sD_3^{d(k)}\\
          &        &        &       &  \cdots\\
          &        &        &       &  \cdots\\
          &        &        &       &  \cdots\\
  0       &0       &0       &\cdots  &^sD_s^s &^sD_s^{s+1} &\cdots&^sD_s^{d(k)}
      }\right],
      \end{eqnarray*}
where
$$
  1\ge{^s}D_j^j\ge \frac{1}{2^{s-j}}, \ 1\le j\le s-1; \ \ {^s}D_s^s =1;$$
and for each $1\le j\le s$,
$$ \max\limits_{1\le l\le d(k)} {^s}D_j^l =1.$$
Moreover,
$$
 \sum\limits_{j=1}^{s}\phi_j(z)\psi_j(\xi)=\sum\limits_{j=1}^{s}{^s\phi_j(z)}\cdot{^s\psi_j(\xi)}.
$$

Let $\tilde \phi_j: ={^s}\phi_j, \tilde \psi_j: ={^s}\psi_j, 1\le j\le s$. Then from the construction, for $1\le j\le s$, $\|\tilde \phi_j\|= 1$ with $\sum\limits_{j=1}^{s}\phi_j(z)\psi_j(\xi)=\sum\limits_{j=1}^{s}{\tilde\phi_j(z)}{\tilde\psi_j(\xi)}.$ Hence
$$
  \|\sum\limits_{j=1}^{s}\phi_j(z)\psi_j(\xi)\|\le \sum\limits_{j=1}^{s}\|\tilde\phi_j\|\|\tilde\psi_j\|\le s\|\tilde \psi\|.
$$
Furthermore, since ${^s}D_j^j\ge \frac{1}{2^{s-j}}$ when $1\le j\le s$,
$$
  \|\sum\limits_{j=1}^{s}\phi_j(z)\psi_j(\xi)\|\ge \max_{1\le j\le s}{^s}D_j^j\cdot\|\tilde\psi_j\|\ge \frac{1}{2^{s-1}}\|\tilde \psi\|.
$$
The proof of Lemma \ref{6} is therefore complete. \qed
\medskip

\noindent{\bf Proof of Lemma \ref{bound}:} Assume by contradiction that there exist families of $\{\phi^\lambda\}$ and $\{\psi^\lambda\}$, such that
 \begin{equation}\label{19}
   \sum\limits_{j=1}^{n-1}\phi^\lambda_j(z)\psi^\lambda_j(\xi)= H^\lambda(z, \xi) + B^\lambda(z, \xi)\langle z, \xi\rangle_\ell
 \end{equation}
 with $\|H^\lambda\|\le C$ while $\|\sum\limits_{j=1}^{n-1}\phi^\lambda_j(z)\psi^\lambda_j(\xi)\|=\lambda\rightarrow \infty$. Applying Lemma \ref{6} to $\phi^\lambda$ and $\psi^\lambda$ if necessary, we can further assume that $\phi^\lambda$ and $\psi^\lambda$ satisfy
$$\|\phi^\lambda\|\le 1, \ \ \ C \|\psi^\lambda\|\le \|\sum\limits_{j=1}^{n-1}\phi_j^\lambda(z)\psi_j^\lambda(\xi)\|=\lambda \le  (n-1)\|\psi^\lambda\|.$$
In special, for each $1\le j\le n-1$,
 \begin{equation}\label{103}\|\phi_j^\lambda\|\le 1, \ \ \frac{1}{n-1}\le \|\frac{\psi_j^\lambda}{\lambda}\|\le \frac{1}{C}.\end{equation}

Dividing both sides of (\ref{19}) by $\lambda$, then one obtains for some polynomial $\tilde B^\lambda$ that
  \begin{equation}\label{20}
   \sum\limits_{j=1}^{n-1}\phi^\lambda_j(z)\frac{\psi^\lambda_j(\xi)}{\lambda}= \frac{H^\lambda(z, \xi)}{\lambda} + \tilde B^\lambda(z, \xi)\langle z, \xi\rangle_\ell.
 \end{equation}
 Since $\phi^\lambda$ and $\psi^\lambda$ satisfy (\ref{103}), we deduce after passing to a subsequence that $\phi^\lambda$ and $\frac{\psi^\lambda}{\lambda}$ converges, say, to polynomials $\phi^\infty$ and $\psi^\infty$. Moreover, the same inequalities in (\ref{103}) pass onto $\phi^\infty$ and $\psi^\infty$ without change, i.e., $\|\phi^\infty\|\le 1$, $\frac{1}{n-1}\le\|\psi^\infty\|\le C$ and $C \|\psi^\infty\|\le \|\sum\limits_{j=1}^{n-1}\phi^\infty_j(z)\psi^\infty_j(\xi)\|\le  (n-1)\|\psi^\infty\|$.

 On the other hand, from (\ref{20}) and Lemma \ref{18}, after passing $\lambda\rightarrow \infty$, there exists some polynomial $B^\infty$ such that
 $$
   \sum\limits_{j=1}^{n-1}\phi^\infty_j(z)\psi^\infty_j(\xi)= B^\infty(z, \xi)\langle z, \xi\rangle_\ell.
 $$
According to Lemma \ref{Hu}, it immediately gives that
$$
  \sum\limits_{j=1}^{n-1}\phi^\infty_j(z)\psi^\infty_j(\xi)=0.
$$
This however contradicts with the fact that $\|\sum\limits_{j=1}^{n-1}\phi^\infty_j(z)\psi^\infty_j(\xi)\|\ge C\|\psi^\infty\|\ge \frac{C}{n-1}$. Therefore, there exists some $\tilde C$ dependent only on $(C, k, m, n)$ such that $ \|\sum\limits_{j=1}^{n-1}\phi_j(z)\psi_j(\xi)\|\\ \le \tilde C$ and hence $\| B\| \le \tilde C$  because of Lemma \ref{18}. \qed
\medskip

With a routine induction process, Lemmas \ref{18} and  \ref{bound}
combined together can be used to show the following:
\medskip

\begin{lemma}\label{5}
 Let $\{\phi_{jr}\}_{j=1}^{n-1}$ and $\{\psi_{jr}\}_{j=1}^{n-1}$ be two families of holomorphic polynomials in $\mathbf C^n$, $1\le r \le m$. Let $H(z, \xi), B(z, \xi)$ be two polynomials in $(z, \xi)$. Suppose that
$$
  \sum\limits_{r=1}^m\bigg(\sum\limits_{j=1}^{n-1}\phi_{jr}(z)\psi_{jr}(\xi)\bigg)\langle z, \xi\rangle_\ell^r= H(z, \xi) + B(z, \xi)\langle z, \xi\rangle_\ell^{m+1}
$$
  and $\|H\|\le C$. Then $\| B\| \le \tilde C$ and
$\|\sum\limits_{j=1}^{n-1}\phi_{jr}(z)\psi_{jr}(\xi)\| \le \tilde C$ for all $1\le r\le m$ with $\tilde
C$ dependent only on $(C, n, m)$ and the degrees of $\phi_{jr},
\psi_{jr}$ for all $1\le r\le m$.
\end{lemma}
\medskip

\section{Proof of the main theorems}

The proof of the main theorems is motivated by the ideas  in
\cite{HZ2} and \cite{EHZ}. Assume $F$ is  not CR transversal to
$M_{\ell}$ at $0$ and $F(U)\not\subset H_{\ell}^{N}$. Assume also $N-n<n-1$ for the moment.

By a result of \cite{BER2}, the set of  points where the
CR transversality holds for  such an $F$ forms an open dense subset
in $M_{\ell}$. Choose a sequence $\{p_j\}\in M_{\ell}$ such that $p_j\rightarrow 0$ and $F$ is CR transversal at each $p_j$ with $ j\ge 1$. Write $q_j: = F(p_j)$. Now for each $j$, applying the normalization
 process to $F$ at $p_j$ as in section 2, we obtain $F^\sharp_{p_j}$ in the following form:
\begin{equation}\label{pj}
      F_{p_j}^\sharp = (f_{p_j}^\sharp, \phi_{p_j}^\sharp, g_{p_j}^\sharp) = ({f_1}_{p_j}^\sharp, \ldots, {f_n}_{p_j}^\sharp, \phi_{p_j}^\sharp, g_{p_j}^\sharp): =T_{p_j}\circ\tau_{F(p_j)}\circ F\circ \sigma_{p_j} ,\ \
      \end{equation}
where
   \begin{eqnarray*}
    &{f}_{p_j}^\sharp(z, w)=z+\frac{i}{2}a^{(1,0)}_{p_j}(z)w +o_{wt}(3),\\
    &\phi_{p_j}^\sharp(z, w)=\phi^{(2,0)}_{p_j}(z) + o_{wt}(2),\\
    &g_{p_j}^\sharp(z, w)=w+o_{wt}(4),
     \end{eqnarray*}
    with the following CR  Gauss-Codazzi equation
    \begin{equation} \label{gauss}
    \langle a^{(1,0)}_{p_j}(z), \bar z\rangle_{\ell} |z|_{\ell}^2 = |\phi^{(2,0)}_{p_j}(z)|^2+ \frac{1}{4}\s S_{p_j}(z).
 \end{equation}
Here $\tau_{F(p_j)}$ is the translation map of $H_\ell^{N}$ sending $F(p_j)$ to 0,
  $\sigma_{p_j}$ is a biholomorphic map   sending 0 to $p_j$ such that $\sigma_{p_j}^{-1}(M_\ell)$ is normalized up to the 4th
  order, and $\s S_{p_j}$ is  the resulting
   Chern-Moser-Weyl curvature function of $M_\ell$ at $p_j$. Note   $\sigma_{p_j}$ depends smoothly on $p_j$.
         Since $F$ is not CR transversal at $0$, $\lim_{j\rightarrow \infty}\lambda_{p_j}=0$ with $\lambda_{p_j}$ defined in (\ref{eqn:022}) for the map $\tau_{F(p_j)}\circ F\circ \sigma_{p_j} $.
       By construction, at each point $p_j$, $F_{p_j}^\sharp$ sends $\sigma_{p_j}^{-1}(M_{\ell})$ into $H_\ell^{N}$.
   We then have for  $(z,u)\approx 0$,
 \begin{eqnarray}\label{eqn}
   &-\Im g_{p_j}^\sharp(z, u+i(|z|_\ell^2+o_{wt}(3))) +|f_{p_j}^\sharp(z, u+i(|z|_\ell^2+o_{wt}(3)))|_\ell^2 + \nonumber\\
   &\   \ \ \ \ \ \ \ \ \ \ \ \ \ \
    + |\phi_{p_j}^\sharp(z, u+i(|z|_\ell^2+o_{wt}(3)))|^2
   =0,
  \end{eqnarray}
 Here $(z, u+i(|z|_\ell^2+o_{wt}(3)))$ is a local parametrization of $\sigma_{p_j}^{-1}(M_{\ell})$ near 0.
  Due to the smooth dependence of $\sigma_{p_j}$ with respect to $p_j$,  the error term $o_{wt}(3)$ depends smoothly on $p_j$.
With an abuse of notation, we shall suppress $\sharp$  and the subindex $j$ of $p$ for the map in (\ref{eqn}).

Given any positive integer $k\ge 2$, collect terms of weighted degree $k$  in the power series expansion of
(\ref{eqn}). We have:
\begin{eqnarray}\label{imp}
  &\Im g_p^{(k)}(z, w) - 2\Re \langle f_p^{(k-1)}(z, w), \bar z\rangle_\ell = (|\phi_{p}(z, w)|^2)^{(k)}\nonumber\\
  +&\ H(g_p^{(r)}|_{0\le r\le k-1}, f_p^{(r)}|_{0\le r\le k-2}, (|\phi_{p}|^2)^{(r)}|_{0\le r\le k-1})
\end{eqnarray}
on $w=u+i|z|_\ell^2$. Here $H$ is a certain bounded polynomial on all its variables. From now on and in what follows,  we use  $C$ in general to represent constants independent of $p$, and use $H(\cdot, \cdot)$ in general to represent polynomials whose norm is bounded by $C$. $C$ and $H$ may be different in different contexts.
\medskip

\begin{lemma}\label{le} Assume that $N-n<n-1$.
For $F_p$ constructed as above and for each $k$, $\|F_p^{(k)}\|\le C$ with $C$ independent of $p$.
\end{lemma}
\medskip

\noindent{\bf Proof of Lemma \ref{le}}:  According to the
normalization procedure conducted in Section 2, $\|g_p^{(k)}\|\le C,
\|f_p^{(k-1)}\|\le C, \|(|\phi_{p}|^2)^{(k)}\|\le C$ automatically
hold when $k\le 4$. Indeed, $\|g_p^{(k)}\|\le 1, \|f_p^{(k-2)}\|\le
1, \|(|\phi_{p}|^2)^{(k-1)}\|=0, k\le 4$ by (\ref{pj}). Moreover,
since $\|\s S_{p}\|\le C$, applying  Lemma \ref{bound}  to
(\ref{gauss}), one has $\|f_p^{(3)}\|\le C,
\|(|\phi_{p}|^2)^{(4)}\|\le C$.

Assuming by induction that $(\|g_p^{(j)}\|, \|f_p^{(j-1)}\|,
\|(|\phi_{p}|^2)^{(j)}\|)$ are all uniformly bounded by some
constant independent of $p$ for $j\le k$,
 we shall show the unform boundedness of $(\|g_p^{(k+1)}\|, \|f_p^{(k)}\|,  \|(|\phi_{p}|^2)^{(k+1)}\|)$. Complexifying (\ref{imp}) at level $k+1$, we obtain
\begin{eqnarray}\label{3}
  &g_p^{(k+1)}(z, w)-\bar g_p^{(k+1)}(\xi, \eta) - 2i \langle f_p^{(k)}(z, w), \xi\rangle_\ell - 2i \langle \bar f_p^{(k)}(\xi, \eta), z\rangle_\ell\nonumber\\
   &= 2i\langle\phi_{p}(z, w), \bar\phi_{p}(\xi, \eta)\rangle^{(k+1)}+ H(z,  \xi, w,\eta)
\end{eqnarray}
which holds on $w-\eta = 2i\langle z, \xi\rangle_\ell$.

Let $L_j=\frac{\p}{\p z_j} + 2i\delta_j\xi_j\frac{\p}{\p w}, 1\le j\le n-1$ with $\delta_j=-1$ when $j\le \ell$ and $\delta_j=1$ with $j\ge\ell+1$. Then $L_j$ is a holomorphic tangent vector field on $w-\eta = 2i\langle z, \xi\rangle_\ell$ for each $j$. Applying $L_j$ onto (\ref{3}), we get
\begin{eqnarray}\label{4}
  &L_jg_p^{(k+1)}(z, w) - 2i \langle L_jf_p^{(k)}(z, w), \xi\rangle_\ell - 2i \bar f_{p, j}^{(k)}(\xi, \eta)\nonumber\\
   &=  2iL_j\langle\phi_{p}(z, w), \bar\phi_{p}(\xi, \eta)\rangle^{(k+1)}+ H(z,  \xi, w, \eta)
\end{eqnarray}
on $w-\eta = 2i\langle z, \xi\rangle_\ell$.

Now we expand $g_p^{(k+1)}, f_p^{(k)}, \langle\phi_{p}(z, w), \bar\phi_p(\xi, \eta)\rangle^{(k+1)}$ in the following manner:
\begin{eqnarray*}
      &g_p^{(k+1)}(z, w)= \sum\limits_{\mu+2\nu=k+1}(g_p)_{\mu\nu}(z)w^\nu;\\
    &f_p^{(k)}(z, w)= \sum\limits_{\mu+2\nu=k}(f_p)_{\mu\nu}(z)w^\nu;\\
    &\langle\phi_{p}(z, w), \bar\phi_p(\xi, \eta)\rangle^{(k+1)}= \sum\limits_{\mu+\gamma+2(\nu+\delta)=k+1}(A_p)_{\mu\gamma\nu\delta}(z, \xi)w^\nu\eta^\delta.
 \end{eqnarray*}
Here $(g_p)_{\mu\nu}$ and $(f_p)_{\mu\nu}$ are homogeneous polynomials of degree $\mu$ in $z$, $(A_p)_{\mu\gamma\nu\delta}$ is a homogeneous polynomial of bi-degree $(\mu, \gamma)$ in $(z, \xi)$.

Let $w=0, \eta = -2i \langle z, \xi\rangle_\ell$ in (\ref{3}). Then
we have
\begin{eqnarray}\label{13}
  &(g_p)_{(k+1)0}(z)-\sum\limits_{\mu+2\nu=k+1}(\bar g_p)_{\mu\nu}(\xi)\eta^\nu - 2i \langle (f_p)_{k0}(z), \xi\rangle_\ell \nonumber\\&- 2i \langle \sum\limits_{\mu+2\nu=k} (\bar f_p)_{\mu\nu}(\xi)\eta^\nu, z\rangle_\ell
   = 2i\sum\limits_{\mu+\gamma+2\delta=k+1}(A_p)_{\mu\gamma0\delta}(z, \xi)\eta^\delta+ H(z,  \xi,\eta)
 \end{eqnarray}
on  $\eta = -2i \langle z, \xi\rangle_\ell$.

Collect terms in $(z, \xi)$ of bi-degree $(k+1, 0)$ and $(k, 1)$  in
(\ref{13}).  By the fact that $\phi_p(0)=\frac{\p\phi_p}{\p
z}(0)=\frac{\p\phi_p}{\p w}(0)=0$ and the definition of
$(A_p)_{\mu\gamma\nu\delta}$,
\begin{eqnarray}\label{133} (A_p)_{k+1,0,0,0}=(A_p)_{k,1,0,0}=(A_p)_{k-1,0,0,1}=0.\end{eqnarray} Then we
have that
$$
  \|(f_p)_{k 0}\|\le C, \ \ \ \|(g_p)_{(k+1) 0}\|\le C.
$$
Hence for each $1\le j\le n-1$,
\begin{eqnarray}\label{8}
&L_jf_p^{(k)}(z, 0) = 2i\delta_j\xi_j\sum\limits_{\mu=k-2}(f_p)_{\mu1}(z)+H(z);\nonumber\\
&L_jg_p^{(k+1)}(z, 0) = 2i\delta_j\xi_j\sum\limits_{\mu=k-1}(g_p)_{\mu1}(z)+H(z).
\end{eqnarray}

Collecting terms in $(z, \xi)$ of bi-degree $(\alpha, \beta), \beta\ge 2$ with $\alpha+\beta =k+1$ in (\ref{13}) gives
\begin{eqnarray}\label{7}
      &-(\bar g_p)_{\beta-\alpha, \alpha}(\xi)\eta^\alpha - 2i\langle z, (\bar f_p)_{\beta-\alpha+1, \alpha-1}(\xi)\eta^{\alpha-1}\rangle_\ell\nonumber\\
    &= 2i \sum\limits_{\theta=0}^{\alpha-2}(A_p)_{\alpha-\theta, \beta-\theta, 0, \theta}(z, \xi)\eta^\theta + H(z,  \xi, \eta)
\end{eqnarray}
with $\eta = -2i \langle z, \xi\rangle_\ell$. Here once again we used the fact that $\phi_p(0)=\frac{\p\phi_p}{\p z}(0)=0$ which implies $(A_p)_{1,(\beta-\alpha-1),0,(\alpha-1)}=(A_p)_{0, \beta-\alpha, 0, \alpha}=0$, so the summation on
 the right hand side runs only till $\theta=\alpha-2$.  Recall from the definition of $A_p$, $(A_p)_{\mu\gamma\nu\delta}(z, \xi) =
 \sum\limits_{j=1}^{N-n} \phi_{p, j}^{(\mu, \nu)}(z, 1)\bar\phi_{p. j}^{(\gamma, \delta)}(\xi, 1)$. Since $N-n< n-1$ by assumption, we immediately have,
  by applying Lemma \ref{5} to (\ref{7}) and by using (\ref{133}), that
\begin{equation}\label{9}
  \|(\bar g_p)_{\beta-\alpha, \alpha}(\xi)\langle z, \xi\rangle_\ell - \langle z, (\bar f_p)_{\beta-\alpha+1, \alpha-1}(\xi)\rangle_\ell \|\le C
\end{equation}
with $\beta\ge 2$, and
 $$
  \| (A_p)_{\mu\gamma0\delta}\|\le C.
$$
Hence from the above inequality,
\begin{equation}\label{10}
  L_j(A_p)(z, \xi, 0, \eta)= 2i\delta_j\xi_j\sum\limits_{\mu+\gamma+2\delta=k-1} (A_p)_{\mu\gamma1\delta}(z, \xi)\eta^\delta +H(z, \xi, \eta).
\end{equation}

Letting $w=0, \eta = -2i\langle z, \xi\rangle_\ell$  and then
substituting (\ref{8}) and (\ref{10}) in (\ref{4}),  we have for
each $1\le j\le n-1$,
 \begin{eqnarray}\label{14}
        &2i\delta_j\xi_j\sum\limits_{\mu=k-1}(g_p)_{\mu1}(z) - 2i\langle 2i\delta_j\xi_j\sum\limits_{\mu=k-2}(f_p)_{\mu1}(z), \xi\rangle_\ell - 2i  \sum\limits_{\mu+2\nu=k+1}(\bar f_{p, j})_{\mu\nu}(\xi) \eta^\nu\nonumber\\
  & = 2i\delta_j\xi_j\sum\limits_{\mu+\gamma+2\delta=k-1} (A_p)_{\mu\gamma1\delta}(z, \xi)\eta^\delta +H(z, \xi, \eta)
    \end{eqnarray}
on $ \eta = -2i\langle z, \xi\rangle_\ell$.  Collect terms in $(z, \xi)$ of bi-degree $(k-1, 1)$ and $(k-2, 2)$  in (\ref{14}). Since $(A_p)_{k-1, 0, 1, 0}=(A_p)_{k-3, 0, 1, 1}=0$, one obtains  that
\begin{eqnarray}\label{12}
  &\|(g_p)_{(k-1)1}\|\le C, \nonumber\\
  &\|2i\delta_j\langle \xi_j (f_p)_{(k-2)1}(z), \xi\rangle_\ell + (\bar f_{p, j})_{(4-k)(k-2)}(\xi)(-2i\langle z, \xi\rangle_\ell)^{k-2}\|\le C.
\end{eqnarray}
Here we have used the convention that $h_{\mu}=0$ if $\mu$ is negative.

Moreover, collecting terms of bi-degree $(\alpha, \beta)$  in $(z, \xi)$ with $\beta\ge 3$ and $\alpha+\beta=k$ in (\ref{14}), one gets for each $1\le j\le n$,
\begin{eqnarray*}
  &(\bar f_{p, j})_{(\beta-\alpha)\alpha}(\xi)(-2i\langle z, \xi\rangle_\ell)^\alpha \\
  &= -\delta_j\xi_j\sum\limits_{\theta=0}^{\alpha-1}(A_p)_{(\alpha-\theta)(\beta-\theta-1) 1\theta}(z, \xi)(-2i\langle z, \xi\rangle_\ell)^\theta +H(z, \xi).
\end{eqnarray*}
Here we use the fact that $(A_p)_{0, \beta-\alpha-1, 1, \alpha}=0$, so the summation on the right hand sides runs only till $\alpha-1$. Applying Lemma \ref{5} onto the above identity as before, we obtain
$  \|(f_p)_{\mu\nu}\|\le C$
for $\mu+2\nu=k$ with $\mu+\nu\ge 3$. When $\mu+2\nu=k\ge 4$ with $\mu+\nu\le 2$, or equivalently, when $\mu=0, \nu=2$, one substitutes the fact that $\|(f_p)_{21}\|\le C$ into (\ref{12}) and gets  $\|(f_p)_{02}\|\le C$ and hence
\begin{equation}\label{11}
   \|(f_p)_{\mu\nu}\|\le C
\end{equation}
for $\mu+2\nu=k$.
Substitute (\ref{11}) into (\ref{9}), then
\begin{equation}\label{15}
  \|(g_p)_{\mu\nu}\|\le C
\end{equation}
for $\mu+2\nu=k+1$ (with $\mu+\nu\ge 3$, which is always fulfilled when $\mu+2\nu=k+1\ge 5$).

Using equation (\ref{3}), we then  have from (\ref{11}) and (\ref{15}),
\begin{equation}\label{1}
 \langle\phi_{p}(z, w), \bar\phi_{p}(\xi, \eta)\rangle^{(k+1)}= H(z,  \xi, w,\eta)
 \end{equation}
 on $w-\eta =2i\langle z, \xi\rangle_\ell$.

 We claim that, for arbitrary $(z,w,\xi,\eta)$, we have
$$
 \langle\phi_{p}(z, w), \bar\phi_{p}(\xi, \eta)\rangle^{(k+1)}= H(z,  \xi, w,\eta).
$$
 Indeed, by (\ref{1}), we have
 \begin{equation}\label{2}
   \sum\limits_{\mu+\gamma+2(\nu+\delta)=k+1}(A_p)_{\mu\gamma\nu\delta}(z, \xi)\bigg(\eta +2i\langle z, \xi\rangle_\ell\bigg)^\nu\eta^\delta = H(z,  \xi, \eta)
 \end{equation}
 near 0. If $\|(A_p)_{\mu\gamma\nu\delta}\|\le C$ does not hold uniformly in $p$, then there exists a smallest integer $\delta_0$
 such that $\|(A_p)_{\mu\gamma\nu\delta_0}\|\rightarrow \infty$ as $p\rightarrow 0$ after passing to a subsequence if necessary.
 Moving the terms with $\delta<\delta_0$ to the right,
  we obtain
$$
    \sum\limits_{\mu+\gamma+2(\nu+\delta_0)=k+1}(A_p)_{\mu\gamma\nu\delta_0}(z, \xi)\big(2i\langle z, \xi\rangle_\ell\big)^\nu = H(z,  \xi).
 $$
Collecting terms in $(z, \xi)$ of bi-degree $(\alpha, \beta)$ with $\alpha+\beta =k+1-2\delta_0$ in the above expression, we get
$$
  \sum\limits_{\theta=0}^\alpha(A_p)_{(\alpha-\theta)(\beta-\theta)\theta\delta_0}(z, \xi)\big(2i\langle z, \xi\rangle_\ell\big)^\theta= H(z,  \xi).
$$
 Recall $(A_p)_{\mu\gamma\nu\delta} = \sum\limits_{j=1}^{N-n} \phi_{p, j}^{(\mu, \nu)}(z, 1)\bar\phi_{p. j}^{(\gamma, \delta)}(\xi, 1)$ by
 definition and $N-n< n-1$. Applying Lemma \ref{5} to the above identity, one deduces $\|(A_p)_{\mu\gamma\nu\delta_0}\|\le C$ for $\mu+\gamma+2(\nu+\delta_0)=k+1$. This is a contradiction! Hence the claim holds.

 The induction is thus complete. Consequently, for each $k$,
$    \|\phi_{p}^{(k)}\|\le C$ with $C$ independent of $p$. We have shown for each fixed $k$,
$\{\|F_{p_j}^{(k)}\|\}_{j=1}^\infty$ is bounded by some constant
independent of $j$.  \qed

\medskip

We are now  in a position to prove Theorem \ref{21} and Theorem \ref{main}, making
use of the result in \cite{MMZ}.
\medskip

\noindent{\bf Proof of Theorem \ref{21}:} If $F$ is CR transversal to $M_\ell$ at $0$, then we are done. Assume $F$ is not CR transversal at $0$.
Then there exists $p_j\rightarrow 0$ such that $F_{p_j}$ as constructed at the beginning of the section satisfies (\ref{pj}). Moreover, for each $k$, $\|F_{p_j}^{(k)}\|\le C$ with $C$ independent of $j$ by Lemma \ref{le}. Following the same trick as in \cite{HZ2}, for each $k$,
$\{F_{p_j}^{(k)}\}_{j=1}^\infty$ converges as $j\rightarrow \infty$ after passing to subsequences,
 to a certain $F^{*(k)}$. By the way these maps were
constructed, the nontrivial formal map $F^*(=(f^*,
\phi^*, g^*)):=\sum\limits_k F^{*(k)}$ sends $M_\ell$ into $H_\ell^{N}$ satisfying
 the following normalization:
\begin{eqnarray*}
    &{f}^*(z, w)=z+o_{wt}(2),\\
    &\phi^*(z, w)= o_{wt}(1),\\
    &g^*(z, w)=w+o_{wt}(4).
  \end{eqnarray*}

According to a result of Meylan-Mir-Zaitsev [MMZ], the formal map
$F^*$ is  convergent. Hence, $F^*$ is a CR immersion
at 0 sending $M_\ell$ into $H_\ell^{N}$. \qed

\medskip

\noindent{\bf Proof of Theorem \ref{main}:} Assume by contradiction that $F$   neither is  CR transversal
 to $M_{\ell}$ at $0$ nor sends $U$ into $H_{\ell}^{N}$. Then there exists a CR immersion $F^*$ sending $M_\ell$ into $H_\ell^{N}$ by Theorem \ref{21}. On the other hand, since any two CR transversal maps between  a Levi-nondegenerate hypersurface and a hyperquadric of the
same signature  differ only by an automorphism of the hyperquadric
(see \cite{EHZ}) when the codimension is less than $\frac{n-1}{2}$,
there exists an automorphism $T$ of $H_\ell^{N}$ such that near
$p_j\approx 0$, and hence at all points in $M_\ell$ near the origin,
$$
   F=T \circ  F^*.
$$
Since $T$ extends to an automorphism  of the projective space
${\mathbf P}^{N}$ and $T(0)=0$, $F$ must be CR transversal at $0$.
  This is a contradiction. \qed
\bigskip

\begin{acknowledgement} Part of work was done when the authors were attending the 7th workshop
on geometric analysis of PDE and several complex variables  at Serra Negra, Brazil. Both authors would like to thank the organizers
for the kind invitation. The second  author is also indebted to Wuhan University for the hospitality during her several visits there.
\end{acknowledgement}


\end{document}